\documentclass[a4paper,10pt]{article}

\usepackage{amssymb,amsmath,amsfonts,amsthm}
\usepackage{pslatex}
\usepackage{a4wide}
\usepackage{cite}

\theoremstyle{definition}

\newtheorem{Definition}{Definition}
\numberwithin{Definition}{section}
\newtheorem{Theorem}[Definition]{Theorem}
\newtheorem{Lemma}[Definition]{Lemma}
\newtheorem{Corollary}[Definition]{Corollary}
\newtheorem{Remark}[Definition]{Remark}

\newcommand{\N}{\mathbb N}
\newcommand{\Q}{\mathbb Q}
\newcommand{\R}{\mathbb R}
\newcommand{\C}{\mathbb C}

\title{Matchings in Balanced Hypergraphs}
\author{Robert Scheidweiler and Eberhard Triesch}
\date{ } 
\begin{document}
 
\maketitle

\begin{abstract}
We give a new proof of K\"onig's theorem and generalize the Gallai-Edmonds decomposition to balanced hypergraphs in two different ways. Based on our decompositions we give two new characterizations of balanced hypergraphs and show some properties of matchings and vertex cover in balanced hypergraphs.
\end{abstract}

\section*{Introduction}

In this article we investigate balanced hypergraphs. Balanced hypergraphs are one possible generalization of bipartite graphs. They were defined by Berge in \cite{berge:1970}. A recent survey about balanced hypergraphs resp. matrices can be found in \cite{cc:2006}. The problem of finding maximum matchings in abitrary hypergraphs is NP-complete, whereas the problem can be solved by linear programming techniques for the class of balanced hypergraphs. But until now there is no polynomial and combinatorial matching algorithm for balanced hypergraphs known.\\ 
The purpose of this article is to deliver a better insight into the matching problem in balanced hypergraphs. In section 1 we define basic notions about hypergraphs, matchings, etc. and list basic results. Moreover we give a new proof of K\"onig's theorem for balanced hypergraphs, which was originally proved by Berge and Las Vergnas \cite{bergelas:1974} and in a more general version by Fulkerson et al. \cite{fuhoffop:1968}. Moreover we show, how our ideas can be used to augment matchings and to estimate the matching number. Section 2 contains our main theorems: two generalizations of the Gallai-Edmonds decomposition. After proving our decompositions we compare them with the classic one (cf. \cite{gallai:1965} and \cite{edmonds:1968}). In the third Section we give two new characterizations of balanced hypergraphs.

\section{Prerequisites}

In this section we define basic notions and start with hypergraphs and graphs.
Let $V=\{v_1,\cdots,v_n\}$ be a finite set and $E=\{e_1,\cdots,e_m\}$ a collection of subsets of $V,$ such that $e\neq\emptyset$ for all $e\in E$ and $\bigcup\limits_{i=1}^m e_i=V.$ The pair $H=(V,E)$ is called hypergraph, the elements $v_i$ of $V$ are the vertices of $H$ and the elements $e_i$ of $E$ are the edges of $H.$ If $|e|\leq 2$ for all $e\in E,$ we call $H$ a graph. For a subset $W\subset V,$ we call the hypergraph $(W,\{e\cap W \mid e\in E, e\cap W \neq \emptyset \})$ the subhypergraph induced by the set $W.$ Furthermore for a subset $F=\{f_1,\cdots, f_k\}\subset E$ we denote  $V_H(F)=V(F)=\bigcup\limits_{i=1}^{k} f_i$ and we call the hypergraph $(\bigcup\limits_{i=1}^{k} f_i,F)$ the partial hypergraph generated by the set $F.$ We define \begin{eqnarray*}H-v&=& (V\setminus{v},\{e\in E \mid v\notin e\}), \\ H\setminus v & = & (V\setminus{v},\{e\setminus\{v\}\mid e\in E, e\setminus\{v\}\neq\emptyset\}) \text{ for all } v\in V \text{ and}\\ H\setminus e & = & \left(\bigcup\limits_{e\in E\setminus \{e\}} e, E\setminus \{e\}\right) \text{ for all } e\in E. \end{eqnarray*} If a hypergraph $\tilde{H}$ is a partial hypergraph of a subhypergraph of $H$ or a subhypergraph of a partial hypergraph of $H$, we say that $\tilde{H}$ is a partial subhypergraph of $H.$ The hypergraph $H^*=(E,\{V_1,\cdots, V_n\})$ with $V_i=\{e\in E\mid v_i\in e\}$ is called the dual hypergraph of $H.$  
Let $\{v_0,v_1,\cdots,v_l\}\subseteq V$ and $\{e_1,\cdots,e_l\}\subseteq E.$
The sequence $P=v_0 e_1 v_1 e_2 \cdots e_l v_l$ is called a path if $v_{i-1},v_i \in e_i$ for $i=1,\cdots, l$ and $v_0,v_1,\cdots v_l$ are pairwise distinct. Moreover we call $l$ the length of the path $P.$ 
The sequence $C=v_0 e_1 v_1 e_2 \cdots e_l v_l$ is called a cycle if $v_{i-1},v_i \in e_i$ for $i=1,\cdots, l,$ $v_0,v_1,\cdots v_{l-1}$ are pairwise distinct and $v_0=v_l.$ In addition we call $l$ the length of the cycle $C.$ The path $P$ resp. the cycle $C$ is called strong, if there is no edge $e_i$ in $P$ resp. $C$ containing three vertices of the path $P$ resp. cycle $C.$ We denote $V(C)=\{v_0,v_1,\cdots,v_{l-1}\}$ and $V(P)=\{v_0,v_1,\cdots,v_l\}.$\\
Now we come to the class of balanced hypergraphs, a generalization of bipartite graphs due to Berge cp. \cite{berge:1970}.
 We call a hypergraph $H$ balanced, if $H$ contains no strong cycle of odd length. The first theorem deals with hereditary properties of balancedness.
\begin{Theorem} \cite{berge:1970} Let $H=(V,E)$ be a balanced hypergraph, then every partial subhypergraph $\tilde{H}$ of $H$ and the dual hypergraph $H^*$ of $H$ are balanced.
\end{Theorem}

Our next step is to define hypergraph edge colorings.
An edge coloring of $H$ in $k$ colors is a function $c:E\rightarrow\{1,\cdots,k\} $ such that $c(e)\neq c(f)$ for all intersecting edges $e,f \in E.$ The sets $C_i=\{e\in E\mid c(e)=i\}$ for $i=1,\cdots, k$ are called color classes. We say that $H$ has the colored edge property if $H$ has an edge coloring in $\Delta(H)$ colors, with $\Delta(H)$ is the maximum degree of a vertex $v\in V.$ \\
The next theorem is again due to Berge. 

\begin{Theorem} \cite{berge:1973} \label{coloredge} Let $H=(V,E)$ be a balanced hypergraph. Then $H$ has the colored edge property.
 \end{Theorem}

\begin{Remark}\label{colorremark}
Balanced hypergraphs have a lot of beautiful coloring properties. It is possible for example 
to color the vertices of balanced hypergraphs in two colors, such that no edge with more than two vertices completely lies in one color class. Berge's proof of theorem \ref{coloredge} uses this property and gives an algorithmic idea how to obtain a minimum edge coloring from proper vertex 2-colorings of balanced hypergraphs. An algorithm to achieve such vertex 2-colorings has been given by Cameron and Edmonds in \cite{caed:1990}. Their algorithm together with Berge's proof yields an efficient procedure to achieve an edge coloring of a balanced hypegraph in $\Delta(H)$ colors (cp. also \cite{cc:2006}). Later we will describe how this procedure can be used to augment matchings in balanced hypergraphs. $\blacksquare$ 
\end{Remark}

Now we turn to matchings and different kinds of optimality criterions for them.
A subset $M\subseteq E$ is called matching of $H$, if the edges $m\in M$ are pairwise disjoint. We say that a matching $M\subseteq E$ is $d$-maximum for a weight function $d:E \rightarrow \N $, if there is no matching $\tilde{M}$ of $H$ with $\sum\limits_{m\in \tilde{M}} d(m) > \sum\limits_{m\in M } d(m)$ and denote the matching number $\curlyvee_d(H)= \sum\limits_{m\in M } d(m).$ If $d\equiv 1,$ we speak of $E$-maximum matchings and denote the matching number by $\curlyvee_E(H)=|M| $ for a maximum matching $M$ with regard to contained edges. Moreover, if $d(e)=|e|$ for all $e\in E,$ we speak of $V$-maximum matchings and denote the matching number by $\curlyvee_V(H)=|V(M)|. $ If we just speak of a maximum matching or matching number, we mean a $V$-maximum matching concerning contained vertices.\\
A subset $S\subseteq V$ is called stable, if $\tilde{S}\nsubseteq e$ holds for every subset $\tilde{S}\subseteq S$ with $|\tilde{S}|\geq 2$ and for all $e\in E.$ Let $d:V\rightarrow \R$ be a weight function. A stable set $S$ is called maximum weight stable set with regard to the weight function $d$, if there is no other stable set $T$ of $H$ with $\sum\limits_{v\in T} d(v) > \sum\limits_{v\in S} d(v).$ If there is only one weight function in consideration, we just speak of a maximum weight stable set.\\
The notions maximum weight stable set and maximum matching are closely related because any $d$-maximum matching $M$ of $H$ is a maximum weight stable set with regard to the weight function $d$ of the dual hypergraph $H^*$ and vice versa.\\  

In order to state K\"onig's theorem for different kind of matching numbers, we have to define different kinds of minimum vertex cover, too.  
Let $x\in  \N^{|V|}.$ Then $x$ is called $d$-vertex cover for a weight function $d:E\rightarrow \N$, if the inequality $$\sum\limits_{v\in e}x_v \geq d(e) $$ holds for every edge $e \in E$. $x$ is called minimum $d$-vertex cover, if there is no vertex cover $\tilde{x}$ with $\sum\limits_{v\in V}x_v >\sum\limits_{v\in e} \tilde{x}_v$ and we denote the $d$-vertex cover number by $\tau_d(H)=\sum\limits_{v\in V}x_v.$ The notions $V$- resp. $E$- vertex covers are defined for the weight function $d(e)=|e|$ resp. $d(e)=1$ for all $e\in E.$
If we just speak of a vertex cover or vertex cover number, we mean a $V$-vertex cover and the $V$-vertex cover number. If the vector $x$ has entries in $\Q$ instead of $\N$, we speak of fractional vertex covers.\\

Now we are ready to state K\"onig's theorem for balanced hypergraphs, which has been proved in parts by Berge and Las Vergnas \cite{bergelas:1974} and Fulkerson et al. \cite{fuhoffop:1968}. We prove it inductively and without the use of linear programming theory. 

\begin{Theorem} \label{vcm} \cite{bergelas:1974}\cite{fuhoffop:1968} Let $H=(V,E)$ be a balanced hypergraph. Then $$\curlyvee_d(H)=\tau_d(H)$$ for all weight functions $d:E\rightarrow \N.$ In particular $\curlyvee_E(H)=\tau_E(H)$ and  $\curlyvee_V(H)=\tau_V(H).$ \end{Theorem}

\proof At first we prove $\curlyvee_E(H)=\tau_E(H).$ Suppose that there is a balanced hypergraph with $\curlyvee_E(H)<\tau_E(H).$ Choose such a hypergraph $H,$ with $|V|+|E|$ minimal. We distinguish two cases:
\begin{enumerate}
 \item There is a $v\in V,$ which is covered by every $E$-maximum matching. Consider $H-v.$ Then it holds:  $\tau_E(H-v)=\curlyvee_E(H-v)=\curlyvee_E(H)-1.$ Now we can construct an $E$-vertex cover of $H$ by taking a minimum $E$-vertex cover $x$ of $H-v$ and setting $x_v=1.$ Then $\sum\limits_{v\in V}x_v=\curlyvee_E(H).$ This is a contradiction.
 \item There is a $E$-maximum matching $M_v$ with $v\notin V(M_v)$ for all $v\in V.$ Choose an abitrary edge $e\in E$ and an $E$- maximum matching $M_v$ with $v\notin V(M_v)$ for every $v\in e.$ Consider the balanced hypergraph $$\tilde{H}=(\bigcup\limits_{v\in e} V(M_v)\cup e, \bigcup\limits_{v\in e}^* M_v \cup\{e\}),$$ in which the union $\bigcup\limits^*$ is a multiset union, i.e., if there are exactly $k$ different matchings $M_{v_1},\cdots,M_{v_k},$ which contain the edge $f,$ the edge $f$ is contained exactly $k$ times in the edge set of $\tilde{H}.$ \\
Since $\Delta(\tilde{H})\leq |e|,$ we can color the edges of $\tilde{H}$ in $|e|$ colors, because balanced hypergraphs have the colored edge property. Let $C_1,\cdots, C_{|e|}$ be the color classes. Note that color classes are also matchings. Then there must be at least one color class, with more than $|M_v|$ edges, because the number of edges of $\tilde{H}$ is $|e||M_v|+1.$ This is a contradiction because the $M_v$ are $E$- maximum matchings. Hence, the situation in case 2 is not possible. 
\end{enumerate}

We have now proved that $\curlyvee_d(H)=\tau_d(H)$ for $d(e)\in \{0,1\}.$ (If $d(e)=0$ for some edges, remove them from the hypergraph and consider the rest.) Now we prove inductively $\curlyvee_d(H)=\tau_d(H)$ and we can assume that $d(e)\geq 1, $ for all $e\in E$ and $d(e)>1$ for at least one $e\in E.$  Suppose that $\curlyvee_d(H)<\tau_d(H)$ for an abritary balanced hypergraph $H$ and choose $\sum\limits_{e\in E}d(e)$ minimal with this property. Define the incidence weight function $p(v):E\rightarrow \N $ for all $v\in V$ with $p(v)(e)=\begin{cases} 0 & v\notin e \\ 1 & v\in e.\end{cases}$ Then it holds $$\curlyvee_{d-p(v)}(H)=\tau_{d-p(v)}(H) \text{ for all }v\in V. $$ We distinguish again two cases:
\begin{enumerate} 

\item There is a vertex $v^*\in V,$ which is contained in a $(d-p(v^*))$-maximum matching. In this case a matching $M^*$ with $\sum\limits_{m\in M^*}(d-p(v^*))(m)+1=\sum\limits_{m\in M^*}d(m)$ exists in $H.$ Furthermore there is a $d$-vertex cover $x^*$ with weight $\tau_{d-p(v^*)}(H)+1.$ (Choose a minimum $(d-p(v^*))$-vertex cover $x$ and set $x_v^*=x_v$ for all $v\in V\setminus\{v^*\}$ and $x^*_{v^*}=x_{v^*}+1.$) Hence both are optimum and $\curlyvee_d(H)=\tau_d(H).$ 

\item For all $v\in V.$ there is no $(d-p(v))$-maximum matching, which covers $v.$ 

Choose an edge $e\in E$ with $d(e)\geq 2.$ There is a $(d-p(v))$-maximum matching $M_v$ with $v\notin V(M_v)$ for every $v\in e.$ Consider the balanced hypergraph $$\tilde{H}=(\bigcup\limits_{v\in e} V(M_v)\cup e, \bigcup\limits_{v\in e}^* M_v \cup\{e\}),$$ in which the union $\bigcup\limits^*$ is again a multiset union. \\
Since $\Delta(\tilde{H})\leq |e|,$ we color the edges of $\tilde{H}$ in $|e|$ colors and let $C_1,\cdots, C_{|e|}$ be the color classes. The sum of all $d$ edge weights of $\tilde{H}$ is at least  $ d(e)+|e| \min\limits_{v\in e} \left\{\sum\limits_{m\in M_v}d(m)\right\}. $
This is the reason why there is a color class $C_i,$ which has a $d$ weight of at least $1+\min\limits_{v\in e}\left\{\sum\limits_{m\in M_v}d(m)\right\}.$ Now choose $v^*\in e$ such that $\sum\limits_{m\in M_{v^*}}d(m)=\min\limits_{v\in e}\left\{\sum\limits_{m\in M_v}d(m)\right\}.$ We can deduce $$\curlyvee_d(H)\geq \sum\limits_{m\in M_{v^*}}d(m)+1 \geq \sum\limits_{m\in M_{v^*}}(d-p(v^*))(m)+1=\tau_{d-p(v^*)}(H)+1\geq \tau_d(H).$$ This achieves the proof. $\blacksquare$
\end{enumerate}

\begin{Remark}\label{matchaug}
The last proof shows, how matchings of balanced hypergraph can be augmented under certain circumstances. If there is an edge $e\in E$ with positive $d$-weight and for every $v\in e$ there is a matching $M_v,$ which does not cover $v,$ we can put all the matchings $M_v,$ for $v\in e,$ and $e$ together in one hypergraph. Then we apply the edge coloring algorithm, which we mentioned above in remark \ref{colorremark}. Then, at least one of the color classes will be a matching of greater $d$-size than the minimum $d$-size of the $M_v,$ $v\in e.$ If all $M_v$ have the same $d$-size, we get a matching of greater $d$-size. $\blacksquare$
\end{Remark}

The colored edge property can also be used to estimate the matching number of a balanced hypergraph. This will be shown in the next theorem.
We have to define the degree $\deg_H(v)=|\{e\in E\mid v\in e\}|$ of a vertex $v$ of a hypergraph $H=(V,E)$ and we denote $\Delta(H)=\max\limits_{v\in V} \{\deg_H(v)\}.$ 

\begin{Theorem}
Let $H=(V,E)$ be a balanced hypergraph and $q\in \N\setminus\{0\}$. If $$\sum\limits_{v\in V}(\Delta(H)-\deg_H(v))\leq q \Delta(H)-1,$$ then $\curlyvee_V(H)\geq |V|-q+1.$
\end{Theorem}
\proof Let $H$ be balanced hypergraph with $\sum\limits_{v\in V}(\Delta(H)-\deg_H(v))\leq q \Delta(H)-1.$ We can color the edges of $H$ in $\Delta(H)$ colors. Suppose that all color classes cover less than $|V|-q+1$ vertices. Then the sum of all vertex degrees of $H$ is at most $(|V|-q)\Delta(H).$ Since $\sum\limits_{v\in V}(\Delta(H)-\deg_H(v))\leq q \Delta(H)-1,$ the sum of all vertex degrees must be at least $\Delta(H)|V|-(q\Delta(H)-1)=\Delta(H)(|V|-q)+1.$ This is a contradiction. $\blacksquare$

An application of K\"onig's theorem yields the following two lemmas, which will be extensively used in the next sections.

\begin{Lemma} \label{matcheq}
Let $H=(V,E)$ be a balanced hypergraph. Let $M$ be a maximum matching of $H$ and let $x$ be a minimum vertex cover of $H.$ Then
$$\sum\limits_{v\in m} x_v =|m| $$ for all $m\in M.$ 
\end{Lemma}
\proof Assume that $\sum\limits_{v\in m^*} x_v >|m^*| $ for an edge $m^*\in M$. The inequality $\sum\limits_{v\in m} x_v \geq |m| $ must hold for the other edges $m\in M\setminus\{m^*\}$ . Otherwise $x$ cannot be a vertex cover. Summing up the weights of all matching edges we get $$\sum\limits_{v\in V} x_v\geq \sum\limits_{m\in M} \sum\limits_{v\in m}x_v> \sum\limits_{m\in M}|m|.$$ This contradicts theorem \ref{vcm}. $\blacksquare$

\begin{Lemma}\label{vc1}
Let $H=(V,E)$ be a balanced hypergraph. Then
$\curlyvee_V(H)-1=\curlyvee_V(H\setminus v)$ if and only if a minimum vertex cover $x$ of $H$ exists with $x_v=1.$
\end{Lemma}
\proof If $\curlyvee_V(H)-1=\curlyvee_V(H\setminus v)$, there is a vertex cover $\tilde{x}$ of $H\setminus v$ with total weight $\curlyvee_V(H)-1$ and one can obtain a vertex cover of $H$ by setting $\tilde{x}_v=1.$ This vertex cover has the same weight as every maximum matching of $H.$ Therefore it is a minimum vertex cover.\\ 
If there is a minimum vertex cover $x$ of $H$ with $x_{v^*}=1$ for a vertex $v^*\in V$, then $\tilde{x}_v=x_v$ for all $v\in V\setminus\{v^*\}$ is a vertex cover of $H\setminus v^*$ (possibly not minimum). Consider a maximum matching $M$ of $H$. \\
Case 1: $v^*\in V(M)$\\ Then we have found a vertex cover and a matching of $H\setminus v^*$ with the same value, namely $\curlyvee_V(H)-1.$ Hence both are optimum and $\curlyvee_V(H)-1=\curlyvee_V(H\setminus v^*).$\\
Case 2: $v^*\notin V(M)$\\ Then we have found a matching with greater weight than a vertex cover in $H\setminus v^*.$ But this is impossible. $\blacksquare$

\section{Decomposition Theorems}
In this section we give two new decomposition theorems for the class of balanced hypergraphs. These theorems generalize the classic Gallai-Edmonds decomposition, which will also be discussed in this section.

\begin{Theorem}\label{GalEd2}
Let $H=(V,E)$ be a balanced hypergraph. We define the sets \begin{eqnarray*}D_H&=&\{v\in V\mid x_v=0 \text{ for all minimum vertex cover } x \text{ of } H\} \\ &=&\{v\in V\mid v \text{ is not covered by every maximum matching } M \text{ of } H\},\\ P_H&=&\{v\in V \mid x_v\geq 2  \text{ for all minimum vertex cover } x \text{ of } H \} \text{ and}\\ M_H&=&V\setminus (P_H\cup D_H).\end{eqnarray*} Then the following conditions hold: \vspace{-2mm}

\begin{enumerate}

\item $\curlyvee_V(H)=\curlyvee_V(H\setminus v)$ for all $v\in D_H.$

\item $\curlyvee_V(H)\leq\curlyvee_V(H\setminus v)$ for all $v\in P_H.$

\item $\curlyvee_V(H)-1=\curlyvee_V(H\setminus v)$ for all $v\in M_H.$

\item There is no edge $e\subseteq D_H.$

\item $|m|\geq 2|m\cap P_H|$ holds for all edges $m$ of maximum matchings of $H.$  \\ $|m|\geq 2|m\cap P_H|+1$ holds for all edges $m$ of maximum matchings of $H$ with $m\cap M_H\neq \emptyset.$

\item Let $v\in D_H.$ Then it holds:\\ 
$\begin{array}{lclclcl}M_{H} &\subseteq& M_{H\setminus v} & & M_{H\setminus v} &\subseteq&  D_H\setminus\{v\}\cup P_H \cup M_H \\
P_H &\subseteq& M_{H\setminus v}\cup  P_{H\setminus v} & & P_{H\setminus v} &\subseteq& D_H\setminus\{v\}\cup P_H \\
 D_H\setminus\{v\}&\subseteq& M_{H\setminus v}\cup  P_{H\setminus v}\cup  D_{H\setminus v}& & D_{H\setminus v} &\subseteq& D_H\setminus\{v\}.\end{array}$ 

\item Let $v\in M_H.$ Then it holds:\\ 
$\begin{array}{lclclcl}M_{H}\setminus\{v\} &\subseteq& M_{H\setminus v}\cup  P_{H\setminus v}\cup  D_{H\setminus v} & & M_{H\setminus v} &\subseteq&  M_H\setminus\{v\} \\
P_H &\subseteq&   P_{H\setminus v} & & P_{H\setminus v} &\subseteq& M_H\setminus\{v\}\cup P_H \\
 D_H&\subseteq& D_{H\setminus v}& & D_{H\setminus v} &\subseteq& M_H\setminus\{v\}\cup D_H.\end{array}$  

\end{enumerate}
\end{Theorem}
\proof
At first we have to show that the two definitions of $D_H$ are equivalent. Let \\$v^*\in \{v\in V\mid v \text{ is not covered by every maximum matching } M \text{ of } H\}$ and assume there is a minimum vertex cover $x$ with $x_{v^*}>0$. Consider a maximum matching $M$ of $H$, which does not contain $v^*$. Then, because of Lemma \ref{matcheq}, we obtain $$\sum\limits_{v\in V}x_v\geq \sum\limits_{m\in M}|m|+x_{v^*}>\sum\limits_{m\in M}|m|.$$ This contradicts theorem \ref{vcm}.\\
Now suppose there is a vertex $v^*$ covered by every maximum matching and $x_{v^*}=0$ for all minimum vertex cover $x$ of $H.$ Then every minimum vertex cover of $H$ is also a vertex cover of $H\setminus v^*$ (possibly not minimum). If $\curlyvee_V(H\setminus v^*) < \curlyvee_V(H)$, there would be a vertex cover with $x_{v^*}=1$ (cf. lemma \ref{vc1} ). Hence $$\curlyvee_V(H\setminus v^*) \geq \curlyvee_V(H).$$  This inequality together with theorem \ref{vcm} implies that any minimum vertex cover $x$ of $H$ must be a minimum vertex cover of $H\setminus v^*$ and $\curlyvee_V(H)=\curlyvee_V(H\setminus v^*).$\\  The inequality \begin{equation}\label{eq1}\sum\limits_{v\in e\setminus \{v^*\}}x_v = \sum\limits_{v\in e}x_v \geq |e|> |e\setminus \{v^*\}| \end{equation} holds for every edge $e\in E$ with $v^*\in e.$ Because of lemma \ref{matcheq} these edges cannot be contained in any maximum matching of $H\setminus v^*$, but two edges must be contained in every maximum matching of $H\setminus v^*$. Otherwise we would get a maximum matching $\tilde{M}$ of $H$, which does not contain $v^*$ or contains more than $\curlyvee_V(H)$ vertices. In both cases we get a contradiction. Therefore no minimum vertex cover $x$ of $H$ exists, which is also a minimum vertex cover of $H\setminus v^*.$ But this is again a contradiction.   
\begin{enumerate}

\item This is implied by the equivalence of the two definitions of $D_H$ and theorem \ref{vcm}.

\item Part two is a direct consequence of lemma \ref{vc1} and the definition of $P_H.$

\item The set $M_H$ is the set of vertices $v$, for which a minimum vertex cover $x$ of $H$ exists with $x_v=1$ (The existance of minimum vertex covers with $x_v=0$ and $x_v=2,$ implies that there is a minimum fractional vertex cover with $x_v=1.$ Then the matching number of $H$ decreases, if we remove $v$ cp. lemma \ref{vc1}. Thus, a vertex cover with $x_v=1$ exists. ) An application of lemma \ref{vc1} yields, that $\curlyvee_V(H)-1=\curlyvee_V(H\setminus v)$ for all $v\in M_H.$ 

\item Suppose there is an edge $e\subseteq D_H$. The edge $e$ cannot be covered by any minimum vertex cover $x$, because $x_v=0$ for all $v\in D_H.$

\item The next part is a direct consequence of lemma \ref{matcheq}, because the equation $\sum\limits_{v\in m} x_v =|m|$ cannot be satisfied for any matching edge $m$ with $|m|<2|m\cap P_H|.$ We can argue analogously for matching edges $m$ with $v\in m\cap M_H.$ Consider a minimum vertex cover of $H$ with $x_v=1.$

\item Consider $H$ and $H\setminus v$ for a vertex $v\in D_H.$ Then any minimum vertex cover $x$ of $H$ is also a minimum vertex cover of $H\setminus v$ (by deleting $x_v=0$). This is the reason why every vertex $v\in M_H$ is contained in $M_{H\setminus v}.$ Moreover there is no vertex in $P_H\cap D_{H\setminus v},$ because such vertices would have weight $0$ in any minimum vertex cover of $H\setminus v$. \\ Suppose that there is a maximum matching $\tilde{M}$ of $H\setminus v,$ which is not a maximum matching of $H.$ There are some (at least two) edges $\tilde{m}_1,\cdots,\tilde{m}_k  \in \tilde{M}$, which contain $v$ in $H.$ This leads to the same contradiction as before (cp. inequality (\ref{eq1})). Therefore \begin{eqnarray*} 
D_{H\setminus v} &=&\bigcup_{\substack{M: M \text{ is maximum} \\ \text{matching of } H\setminus v}} \left( (V\setminus \{v\})\setminus V(M)\right)\\
 & \subseteq & \bigcup_{\substack{M: M \text{ is maximum} \\ \text{matching of } H}} \left( V\setminus V(M)\right)\setminus\{v\}
\\ & = & D_H\setminus\{ v\}. \end{eqnarray*}

\item Now we turn to vertices $v\in M_H$ and consider $H$ and $H\setminus v.$ Then any minimum vertex cover $x$ of $H\setminus v$ is also a minimum vertex cover of $H$ (by setting $x_v=1$) and maximum matchings of $H$ are maximum matchings of $H\setminus v$ because of part three. Therefore any vertex $v\in P_H$ is contained in $P_{H\setminus v}$ and any vertex $v\in D_H$ is contained in $D_{H\setminus v}.$ $\blacksquare$

\end{enumerate}

The classic Gallai-Edmonds decomposition divides the set of vertices of a graph $H$ into these three sets:
\begin{eqnarray*}D_H&=&\{v\in V\mid v \text{ is not covered by every maximum matching } M \text{ of } H\}, \\ A_H&=&\{v\in V\setminus D_H \mid v \text{ is a neighbor of a vertex } \tilde{v}\in D_H \}\text{ and} \\C_H&=&V\setminus (A_H\cup D_H).\end{eqnarray*}
We have found out, that only part 1, part 4 and a result similar to part 6 of theorem \ref{GalEd2} are valid, if we use this classic decomposition for balanced hypergraphs.
Nevertheless the two different decompositions $D_H, A_H, C_H$ and $D_H,P_H,M_H$ are equal for example in the case of bipartite graphs and in the class of factor critical balanced hypergraphs (i.e. the hypergraph $H\setminus v $ has a perfect matching for all $v\in V$). Which is the greatest subclass of balanced hypergraphs, for which these decompositions are equal? An answer to this question would be very interesting.

\begin{Theorem}\label{seteq}
Let $H=(V,E)$ be a balanced hypergraph and let the sets $D_H, A_H, P_H, C_H$ and $M_H$ be defined as before. Then $A_H=P_H$ and $C_H=M_H$ if:  \vspace{-2mm}

\begin{enumerate}
\item $M_H=\emptyset,$ in particular if $H$ is factor critical. 

\item $H$ is a bipartite graph.

\end{enumerate}
\end{Theorem}

\proof

\begin{enumerate}

\item If $H$ is factor critical the equation $\curlyvee_V(H)\leq\curlyvee_V(H\setminus v)$ holds for all $v\in V.$ Therefore $v\in P_H$ or $v\in D_H.$ Hence $M_H=\emptyset.$ \\ 
Now we consider an abitrary balanced hypergraph with $M_H=\emptyset$ and suppose that there is vertex in $C_H\cap P_H.$ This vertex would not have any neighbors in $D_H$ and we would obtain an edge $m$ of a maximum matching with $ \sum\limits_{v\in m}x_v > |m|.$ This contradicts lemma \ref{matcheq}. Hence we conclude $$ C_H=M_H=\emptyset\text{ and } A_H = P_H. $$

\item $H$ is a bipartite graph. By means of its classic Gallai Edmonds decomposition (cf. \cite{gallai:1965} and \cite{edmonds:1968} ) we can construct a vertex cover of the graph $H.$ Set $x_v=1$ for all $v\in C_H,$  $x_v=2$ for all $v\in A_H$ and $x_v=0$ for all $v\in D_H.$ This is a vertex cover of $H,$ because edges in graphs have at most cardinality $2$ and there is no edge in the subgraph induced by the vertices of $D_H.$ Furthermore $x$ has the same weight as any maximum matching of $H.$ Therefore it is a minimum vertex cover. Summarizing we deduce $C_H \subseteq M_H.$ Now consider an abitrary vertex $v\in A_H.$ If we set $x_v=1,$ there would be an edge between $v$ and $D_H,$ which could not be covered by $x$, because vertices of $D_H$ always have weight $0.$ Hence, we obtain $A_H\subseteq P_H.$ This completes the proof.$\blacksquare$\\   

\end{enumerate}

In our next step we study $E$-maximum matchings and $E$-vertex cover.

\begin{Theorem}\label{GalEd1}
Let $H=(V,E)$ be a balanced hypergraph. We define the sets \begin{eqnarray*}F_H&=&\{v\in V\mid x_v=0 \text{ for all minimum $E$-vertex cover } x \text{ of } H\} \\ &=&\{v\in V\mid v \text{ is not covered by every $E$-maximum matching } M \text{ of } H\},\\ Q_H&=&\{v\in V \mid x_v=1  \text{ for all minimum $E$-vertex cover } x \text{ of } H \} \text{ and}\\ N_H&=&V\setminus (Q_H\cup F_H).\end{eqnarray*} Then the following conditions hold: \vspace{-2mm}

\begin{enumerate}

\item $\curlyvee_E(H)=\curlyvee_E(H\setminus v)$ for all $v\in F_H.$

\item $\curlyvee_E(H)<\curlyvee_E(H\setminus v)$ for all $v\in Q_H,$ which do not lie in an edge of size one. 

\item $\curlyvee_E(H)=\curlyvee_E(H\setminus v)$ for all $v\in N_H.$

\item There is no edge $e\subseteq F_H.$

\item There is no edge $m\in M$ with $m\subseteq N_H\cup Q_H $ and $m\cap Q_H\neq\emptyset$ for every $E$-maximum matching $M.$

\item Let $v\in F_H.$ Then it holds:\\ 
$\begin{array}{lcl}N_{H} &=& N_{H\setminus v}   \\
Q_H &=& Q_{H\setminus v}\\
 F_H\setminus\{v\}&=& F_{H\setminus v}.\end{array}$ 

\item Let $v\in N_H.$ Then it holds:\\ 
$\begin{array}{lclclcl}N_{H}\setminus\{v\} &\subseteq& N_{H\setminus v}\cup  Q_{H\setminus v}\cup  F_{H\setminus v} & & N_{H\setminus v} &\subseteq&  N_H\setminus\{v\} \\
Q_H &\subseteq&   Q_{H\setminus v} & & Q_{H\setminus v} &\subseteq& N_H\setminus\{v\}\cup Q_H \\
 F_H&\subseteq& F_{H\setminus v}& & F_{H\setminus v} &\subseteq& N_H\setminus\{v\}\cup F_H.\end{array}$  

\end{enumerate}
\end{Theorem}
\proof At first we show that the two definitions of $F_H$ are equivalent. It is obvious that \\ $\{v\in V\mid x_v=0 \text{ for all minimum $E$-vertex cover } x \text{ of } H\} \supseteq \{v\in V\mid v \text{ is not covered by every $E$-maximum matching }\},$ because of K\"onig's theorem. Suppose now that there is a vertex $v$ covered by every $E$-maximum matching and $x_v=0$ for all minimum $E$-vertex cover. Then we consider $H-v$ and construct such as in the proof of theorem \ref{vcm} (case 1 of the induction basis) a vertex cover $x$ of $H$ with $x_v=1,$ which yields a contradiction. 

\begin{enumerate}

\item Let $v\in F_H,$ take an $E$- maximum matching $M$, which does not cover $v,$ and a minimum $E$-vertex cover $x$ of $H.$ They are also $E$-vertex cover resp. $E$- maximum matching in $H\setminus v.$ Therefore $\curlyvee_E(H)=\curlyvee_E(H\setminus v).$

\item Let $v\in Q_H.$ Suppose that $\curlyvee_E(H)\geq\curlyvee_E(H\setminus v).$ Then we obtain an $E$-vertex cover of $H\setminus v$ with weight $\curlyvee_E(H).$ This $E$-vertex cover cannot be an $E$-vertex cover of $H$ (by setting $x_v=0$), because $v\in Q_H.$ Therefore $v$ lies in an edge of size one. 

\item Let $v\in N_H.$ take an $E$- maximum matching $M$ and a minimum $E$-vertex cover $x$ of $H$ with $x_v=0.$  They are also minimum $E$-vertex cover resp. $E$- maximum matching in $H\setminus v,$ because $v$ does not lie in an edge of size one. Therefore $\curlyvee_E(H)=\curlyvee_E(H\setminus v).$

\item Cp. theorem \ref{GalEd2}.

\item This assertion is implied by K\"onig's theorem \ref{vcm}.

\item Let $v\in F_H.$ $E$- maximum matchings of $H$ and a minimum $E$-vertex cover of $H$ are also $E$-vertex cover resp. $E$- maximum matching of $H\setminus v.$ Therefore $N_{H} \subseteq N_{H\setminus v}$ and $F_H\setminus\{v\}\subseteq F_{H\setminus v}.$ Moreover there cannot be an $u\in Q_H\setminus Q_{H\setminus v},$ because then there would be an $E$-vertex cover $x$ of $H\setminus v$ with $x_u=0$ and we would obtain an $E$-vertex cover $x$ of $H$ with $x_u=0,$ which is not possible. Hence, $Q_{H} \subseteq Q_{H\setminus v}.$ Summarizing we obtain the three equalities.

\item Let $v\in N_H.$ Then $E$-maximum matchings of $H$ are $E$-maximum matching of $H\setminus v$ and minimum $E$-vertex cover of $H\setminus v$ are minimum $E$-vertex cover of $H.$ Therefore $N_{H}\setminus\{v\} \supseteq N_{H\setminus v}$ and $Q_H \subseteq   Q_{H\setminus v}.$ Moreover $F_H\subseteq F_{H\setminus v}.$ Let $u\in Q_{H\setminus v}, $ then $x_u=1$ for all minimum $E$-vertex cover $x$ of $H\setminus v.$ This is the reason, why $u\in N_H\setminus\{v\}\cup Q_H.$ Let $w\in F_{H\setminus v}, $ then $x_w=0$ for all minimum $E$-vertex cover $x$ of $H\setminus v.$ This is the reason, why $w\in N_H\setminus\{v\}\cup F_H.$ $\blacksquare$
\end{enumerate}

\begin{Theorem}\label{seteq2}
Let $H=(V,E)$ be a balanced hypergraph and let the sets $D_H, F_H, A_H, Q_H, C_H$ and $N_H$ be defined as before. Then $D_H=F_H,$ $A_H=Q_H$ and $C_H=N_H$ if $H$ is a bipartite graph.

\end{Theorem}

\proof It is obvious that $D_H=F_H.$ Again by means of the classic Gallai Edmonds decomposition (cf. \cite{gallai:1965} and \cite{edmonds:1968} ) we can construct an $E$-vertex cover of the graph $H.$ Set $x_v=1$ for one class of the bipartition of $C_H,$  $x_v=1$ for all $v\in A_H$ and $x_v=0$ for all $v\in F_H.$ This is again an $E$-vertex cover of $H,$ because edges in graphs have at most cardinality $2$ and there is no edge in the subgraph induced by the vertices of $F_H.$ Furthermore $x$ has the same weight as any $E$-maximum matching of $H.$ Therefore it is a minimum $E$-vertex cover. Summarizing we deduce $C_H \subseteq N_H.$ Now consider an abitrary vertex $v\in A_H.$ If we set $x_v=0,$ there would be an edge between $v$ and $F_H,$ which could not be covered by $x$, because vertices of $F_H$ always have weight $0.$ Hence, we obtain $A_H\subseteq Q_H.$ This completes the proof.$\blacksquare$

\section{Applications}

In this section we give some different applications of our decomposition, including two new characterizations of balanced hypergraphs similar to Berge's characterizations (cp. \cite{berge:1989}).

\begin{Corollary}\label{charac}
The following statements are equivalent:
\begin{enumerate}
 \item $H=(V,E)$ is balanced.
\item There is no edge $e\subseteq D_{\tilde{H}}$ for every partial subhypergraph $\tilde{H}$ of $H$.

\end{enumerate}
\end{Corollary}
\proof In theorem \ref{GalEd2} we have just seen that for balanced hypergraph there is no egde $e\subseteq D_H$, but we give a second proof here, which is similar to the proof of theorem \ref{vcm}:\\
Suppose there is an edge $e\subseteq D_H.$ There is a maximum matching $M_v$ with $v\notin V(M_v)$ for every $v\in e.$ Consider the balanced hypergraph $$\tilde{H}=(\bigcup\limits_{v\in e} V(M_v)\cup e, \bigcup\limits_{v\in e}^* M_v \cup\{e\}),$$ in which the union $\bigcup\limits^*$ is again a multiset union. \\
Since $\Delta(\tilde{H})\leq |e|,$ we color the edges of $\tilde{H}$ in $|e|$ colors and let $C_1,\cdots, C_{|e|}$ be the color classes. The sum of all vertex degrees of $\tilde{H}$ is $$|e|+|e|\cdot |V(M_v)|.$$ Suppose that the strict inequality $|V(C_i)|<|V(M_v)|$ holds for a color class $C_i.$ Then there must be another color class $C_j$ with $|V(C_j)|>|V(M_v)|,$ because $\sum\limits_{i=1}^{|e|}|V(C_i)|=|e|+|e|\cdot |V(M_v)|.$ This is a contradiction, since any color class is a matching and cannot cover more vertices than a maximum matching.\\
Hence every color class $C_i, i=1,\cdots,|e|$ contains exactly $|V(M_v)|$ vertices. This is again a contradiction because then $\sum\limits_{i=1}^{|e|}|V(C_i)|=|e|\cdot |V(M_v)|\neq |e|+|e|\cdot |V(M_v)|.$\\

For the other direction suppose that $H$ is not balanced. There is a partial subhypergraph $C=(\tilde{V},\tilde{E})$ of $H,$ which is a strong odd cycle. Hence $|\tilde{V}|=2k+1,$ for a $k\in \N\setminus\{0\}$ and $|e|=2$ for all $e\in \tilde{E}.$ Furthermore every vertex of $C$ is contained in $D_C,$ so we get at least three edges contained in $D_C.$ $\blacksquare$

\begin{Remark}\label{DHweight}
Part four of theorem \ref{GalEd2} is also valid, if we speak of maximum matchings with regard to any weight function $d:E\rightarrow \N\setminus\{0\}.$ We can deduce that there is no edge $$e\subseteq\{v\in V\mid v \text{ is not covered by every $d$-maximum matching } M \text{ of } H\}.$$ The proof is the same as in corollary \ref{charac} (compare also the proof of theorem \ref{vcm} and remark \ref{matchaug}), replace $|V(\cdots)|$ by $\sum\limits_{e\in \cdots}d(e).$ $\blacksquare$
\end{Remark}

\begin{Corollary}\label{charac2}
The following statements are equivalent:
\begin{enumerate}
 \item $H=(V,E)$ is balanced.
\item Any vertex $v\in V$ is contained in an edge $e,$ for which $S\cap e\neq\emptyset$ holds for every maximum weight stable set $S$ of $H$ with regard to a weight function $d:V\rightarrow \N\setminus\{0\}.$
\end{enumerate}

\end{Corollary}

\proof We consider part four of theorem \ref{GalEd2}, remark \ref{DHweight} and their consequences in the dual hypergraph $H^*.$ Matchings become stable sets and the set $D_H$ becomes the set of all edges not being covered by every maximum weight stable set, say $D^*_H.$ Hence, there is no vertex solely contained in edges of $D^*_H.$\\
  For the other direction suppose that $H$ is not balanced. There is a partial subhypergraph $C=(\tilde{V},\tilde{E})$ of $H,$ which is a strong odd cycle. Hence $|\tilde{V}|=2k+1,$ for a $k\in \N\setminus\{0\}$ and $|e|=2$ for all $e\in \tilde{E}.$ Set $d(v)=1$ for all $v\in V.$ A maximum weight stable set $S$ with regard to $d,$ has weight (and cardinality) $k.$ Moreover there is an edge $e\in E,$ which is not covered by $S.$ This implies, that there is no edge covered by every maximum weight stable set. $\blacksquare$


\begin{thebibliography}{CCV06}

\bibitem[Ber70]{berge:1970}
Claude Berge.
\newblock Sur certains hypergraphes g\'{e}n\'{e}ralisant les graphes
  bipartites.
\newblock {\em Combinatorial Theory and its Applications I, Colloquium of
  Mathematical Society Janos Bolyai}, 4:119--133, 1970.

\bibitem[Ber73]{berge:1973}
Claude Berge.
\newblock Notes sur les bonnes colorations d'un hypergraphe. colloque sur la
  th{\'{e}}orie des graphes.
\newblock {\em Cahiers Centre \'{E}tudes Rech. Op\'{e}r.}, 15:219--223, 1973.

\bibitem[Ber89]{berge:1989}
Claude Berge.
\newblock {\em Hypergraphs - Combinatorics of Finite Sets}.
\newblock North-Holland, 1989.

\bibitem[BV70]{bergelas:1974}
Claude Berge and Michel~Las Vergnas.
\newblock Sur un th\'{e}or\`{e}me du type k\"onig pour hypergraphes.
\newblock {\em Annals of the New York Academy of Science}, 175:32--40, 1970.

\bibitem[CCV06]{cc:2006}
Michele Conforti, G\'{e}rard Cornu\'{e}jols, and Kristina Vu\v{s}kovi\'{c}.
\newblock Balanced matrices.
\newblock {\em Discrete Mathematics}, 306:2411--2437, 2006.

\bibitem[CE90]{caed:1990}
Kathie Cameron and Jack~R. Edmonds.
\newblock Existentially polytime theorems.
\newblock {\em Dimacs Series in Discrete Mathematics and Theoretical Computer
  Science}, 1:83--100, 1990.

\bibitem[Edm68]{edmonds:1968}
Jack~R. Edmonds.
\newblock Maximum matching and a polyhedron with 0,1-vertices.
\newblock {\em J. Res. Nat. Bur. Standards Sect.}, B:125--130, 1968.

\bibitem[FHO74]{fuhoffop:1968}
Delbert~R. Fulkerson, Alan~J. Hoffman, and Rosa Oppenheim.
\newblock On balanced matrices.
\newblock {\em Mathematical Programming Study}, 1:120--132, 1974.

\bibitem[Gal65]{gallai:1965}
Tibor Gallai.
\newblock Maximale systeme unabh{\"a}ngiger kanten.
\newblock {\em Magyar Tud. Akad. Mat. Kutat\'{o} Int. K\"ozl.}, 9:401--413,
  1965.

\end{thebibliography}
\end{document}